\documentclass[11pt,reqno]{amsart}
\usepackage{amsmath, latexsym, amsfonts, amssymb,amsthm, amscd,epsfig,enumerate, color}

\usepackage{tikz}

\advance\hoffset -.75cm

\usepackage{amsfonts}
\usepackage{latexsym}
\usepackage{amsmath}
\usepackage{amssymb}

\newcommand{\R}{\mathbb R}
\newcommand{\Z}{\mathbb Z}
\newcommand{\N}{\mathbb N}
\newcommand{\Prob}{\mathbb P}

\newcommand{\Rd}{{\mathbb R}^d}
\newcommand{\Zd}{{{\mathbb Z}^d}}

\newtheorem{teo}{Theorem}[section]

\newtheorem{rem}[teo]{Remark}

\newtheorem{defn}[teo]{Definition}

\begin{document}

\title[]
{Combination versus sequential monotherapy in chronic HBV infection: a mathematical approach}

\author[D.~Bertacchi]{Daniela Bertacchi}
\address{D.~Bertacchi,  Universit\`a di Milano--Bicocca
Dipartimento di Matematica e Applicazioni,
Via Cozzi 53, 20125 Milano, Italy
}
\email{daniela.bertacchi\@@unimib.it}

\author[F.~Zucca]{Fabio Zucca}
\address{F.~Zucca, Dipartimento di Matematica,
Politecnico di Milano,
Piazza Leonardo da Vinci 32, 20133 Milano, Italy.}
\email{fabio.zucca\@@polimi.it}

\author[S.~Foresti]{Sergio Foresti}
\address{S.~Foresti, Division of Infectious Diseases
Department of Internal Medicine
"San Gerardo" Hospital, 
20900 Monza, Italy}
\email{s.foresti\@@hsgerardo.org}

\author[D.~Mangioni]{Davide Mangioni}
\address{D.~Mangioni, Division of Infectious Diseases
Department of Internal Medicine
"San Gerardo" Hospital, Universit\`a di Milano-Bicocca
20900 Monza, Italy}
\email{d.mangioni\@@campus.unimib.it}

\author[A.~Gori]{Andrea Gori}
\address{A.~Gori, Division of Infectious Diseases
Department of Internal Medicine
"San Gerardo" Hospital, Universit\`a di Milano-Bicocca
20900 Monza, Italy}
\email{a.gori\@@unimib.it}

\date{\today}

\begin{abstract}
Sequential monotherapy is the most widely used therapeutic approach in the treatment of
HBV chronic infection. Unfortunately, under therapy, in some patients the hepatitis virus mutates and 
gives rise to variants which are drug resistant. 
We conjecture that combination therapy is able to delay drug resistance for a longer time than
sequential monotherapy. To study the action of these two therapeutic approaches in the event
of unknown mutations and to explain the emergence of drug resistance,
we propose a stochastic model for the infection within a patient which is treated with two drugs,
either sequentially or contemporaneously, and develops a two-step mutation which is resistant to
both drugs. We study the deterministic approximation of our stochastic model and give a biological
interpretation of its asymptotic behaviour.
We compare the time when this new strain first reaches detectability in the serum viral load.
Our results show that the best choice is to start an early combination therapy, which allows to
stay drug-resistance free for a longer time.
\end{abstract}
\maketitle

\baselineskip .6 cm



\section{Introduction}

Hepatitis B virus is the cause of one of the most common infections in the world, with around 400 
million people infected according to the WHO assessment \cite{cf:I}. Its consequences, such as the liver 
cirrhosis and the hepatocellular carcinoma, affecting approximately 25\% of patients with chronic 
HBV infection, are important causes of morbidity and mortality worldwide \cite{cf:I, cf:II,cf:1,cf:3}.

The progression of the liver disease relies on environmental factors (co-infection with other viruses 
such as HCV or HIV, alcohol abuse), host factors (age, immune status), and viral factors. As to the
latter, a great importance is assumed by the viral load and the presence of HBV genomic mutations, that
can significantly affect the response to the treatment \cite{cf:III}.
Antiviral therapy has demonstrated to have an important role in reducing both the markers of liver 
disease progression (with new generation nucleoside analogues, viral suppression, i.e.~the
undetectability of HBV DNA, can be reached in about 95\% of  all the cases \cite{cf:IV}) and, 
above all, the evolution to cirrhosis and hepatocellular carcinoma \cite{cf:3}.
 It has been shown that high serum levels of 
HBV-DNA ($> 10^4$ copies/ml) are predictive of an increased risk of hepatocellular carcinoma, 
regardless of the transaminase level and the presence of liver cirrhosis \cite{cf:2}: 
antiviral therapy should therefore be considered a preventive and anti-cancer therapy \cite{cf:3}.

Over recent years several antiviral drugs specific for the virus enzyme DNA polymerase have entered 
in the clinical practice, but so far they are mainly used in a sequential monotherapy
(that is, only one drug is employed at a time and if an increase in viral load is
observed, then the drug is substituted with another one). 
Moreover, interferon monotherapy is still a mainstay for mild to moderate chronic infections \cite{cf:VIII}. 
Some unsolved problems remain today:
the clearance, or viral eradication (certified by the disappearance of HBsAg or cccDNA), which would 
allow the discontinuation of the therapy, is obtained in only 10\% of treated patients \cite{cf:VI}; 
to improve this result, different steps of HBV replicative cycle are being considered as possible 
target for the development of new drugs \cite{cf:VII}.
On the other hand, the long term treatment
with sequential monotherapy \cite{cf:VIII, cf:IX, cf:X}, is likely to select drug resistant strains,
because of the onset of viral breakthrough and mutant escape (that is the emergence of variant viruses)
under the selective pressure of medications \cite{cf:XI, cf:XII}.

As stated in the international guidelines, the combination therapy 
(that is, therapy with a cocktail of drugs)
has a well-defined role in the
treatment of chronic hepatitis B, especially as salvage therapy and first line application in selected
categories of patients, such as those with decompensated cirrhosis, HIV co-infection, pre-existing
mutations and post-liver transplantation \cite{cf:VIII, cf:IX, cf:X, cf:XIII}. 
The guidelines, however, are based only on published results, and to date, only a few studies
have attempted to compare the 
effectiveness of the two different approaches (combination versus monotherapy) in treatment-naive
patients. These studies did show promising results in favor of combination therapy 
(\cite{cf:XIV, cf:XV, cf:XVI, cf:XVII}).

We conjecture that combination therapy, as a  therapy for chronic HBV infection,
is preferable to sequential monotherapy.
Our conjecture is based on the following observation: 
despite being a stable DNA virus, hepatitis B virus replicates in the host cell through a 
RNA phase \cite{cf:5}, 
that makes it similar to the Human Immunodeficiency Virus (HIV). 
The natural history of HIV infection has been radically transformed by the introduction of 
combination antiviral therapy (Highly Active Antiretroviral Therapy, HAART) \cite{cf:6},
 therefore we
believe that this kind of therapy should be beneficial also for HBV treatment.
In particular we infer that under combination therapy the appearance of drug resistance
is less likely, or at least, occurs after a longer time.
We believe it is important to compare the outcome of combination therapy and of
sequential monotherapy in the treatment of chronic HBV infection using a mathematical model,
which takes into account the emergence of mutants of the virus.

During the last two decades and more, mathematical models have been proposed to describe
viral infections (such as HBC, HCV, HIV) both from a single patient's point of view, which is
our point of view in this paper, and from an epidemics point of view.
The standard mathematical model for chronic HBV infection within a single patient
is given by a three ordinary differential equations system (in short, ODEs system):
\begin{equation}\label{eq:stdmodel}
 \begin{cases}
   \dot T=\lambda-\delta^\prime T - \alpha VT; &\\
   \dot Y=\alpha VT-\delta Y;&\\
   \dot V=pY - c V;&\\
  \end{cases}
\end{equation}
where $T$ denotes target (uninfected) hepatocytes which are produced at rate $\lambda$,
die at rate $\delta^\prime T$ and are infected at a rate $\alpha VT$; $Y$ denotes infected
hepatocytes which die at rate $\delta Y$ ($\delta$ possibly larger than $\delta^\prime$) and $V$
denotes free virus. Virions are produced from infected cells at rate $pY$ and are cleared 
from bloodstream at rate $cV$. The system has no exact solutions but approximate solutions
may be obtained under different hypotheses.

If the patient undergoes treatment, the drug acts against reproductions and/or infections.
In the model this is explained by introducing the efficacy parameters $\varepsilon\in[0,1]$
(\cite{cf:Nowak96,cf:TsiangRooney}) and $\eta\in [0,1]$ (\cite{cf:ADEFOVIR-LAMIVUDINA, cf:Lewin})
and replacing
%
%
the previous system by
\begin{equation}\label{eq:basicmodel}
 \begin{cases}
   \dot T=\lambda-\delta^\prime T - (1-\eta)\alpha VT;&\\
   \dot Y=(1-\eta)\alpha VT-\delta Y;&\\
   \dot V=(1-\varepsilon)pY - c V.&\\
  \end{cases}
\end{equation}
The parameters $\varepsilon$  and $\eta$ respectively
indicate the ability of treatment to block viral production or infections
  (the extremal values $0$ and $1$ refer  to no ability
at all and complete ability, respectively). Of course the behaviour under no treatment is
retrieved setting $\eta=\varepsilon=0$.

This model has been used in several papers to fit experimental data under
treatment with different drugs, leading to estimates for some of the 
parameters involved, mainly for $\varepsilon$, $c$ and $\delta$ 
(\cite{cf:ADEFOVIR-LAMIVUDINA, cf:Entecavir, cf:Lamivudina, cf:lau-tsiang, cf:Lewin,
cf:TsiangRooney, cf:Neumann}).
It has also been modified to model HIV infection (\cite{cf:HIV}) or acute
HBV infection (\cite{cf:CRNP}).

To our knowledge, none of the existing models gives an explanation of the emergence
of drug-resistant strains of the virus, which are usually not observed without therapy. 
We want to modify the previous model, according to some biological observations.
Within a single infected person, genetically different viral variants can coexist: that is because of 
the same nature of HBV, which reproduces itself through a pre-genomic RNA phase, accounting for the 
presence of heterogeneous strains in the viral progeny \cite{cf:5, cf:XII}. 
The mutations are pre-adaptive, linked 
to random errors of the polymerase enzyme, especially in the RNA phase of the replicative process. 
Hence, the antiviral therapy does not determine any resistant strains, but simply selects them 
\cite{cf:5, cf:XII, cf:8, cf:9}.
Moreover, in treatment-naive patients a minority of drug-resistant species can also be found. 
These variants could be the result of the persistence of drug-resistant viruses transmitted by a 
treated individual, or may have been generated ex novo during the course of the original infection
\cite{cf:XII}. 

In order to keep the model simple and to focus on the first mutations which appear under
sequential monotherapy, we think of the example of a
treatment-naive patient which undergoes sequential monotherapy, with two drugs at hand 
(drug A and drug B)
and who will develop drug-resistance to both. 
Before therapy, he/she has a high serum viral load of HBV virions, of a type
(usually the wild-type) that we call \textit{type 1}. During therapy with drug A,
the viral load decreases (eventually reaching undetectability) but in the long run a drug-resistant
variant emerges, namely we observe the rise of the serum viral load of what we denote by
\textit{type 2} virions. This forces us to change therapy and resort to drug B, which is able
to strike type 2 down, but eventually another mutation arises (\textit{type 3} mutants), 
which escapes the action of drug B and serum viral load starts increasing again.
We identify the moment when we detect type 3 as the moment of therapeutic failure
(with the current tests at hand, this is usually when serum viral load reaches 20 copies/ml). Of course
we may have at our disposal another drug capable of fighting this mutation too, and the process
could continue, but we are interested in those patients which may eventually develop drug-resistance
and a two-step mutation is the simplest case where one can test the difference between combination
therapy and sequential monotherapy.
Our main question is whether providing drug A and drug B together can significantly delay the
instant of therapeutic failure (if so, it would be an argument in favour of combination therapy).

In Section~\ref{sec:stochmodel} we describe a stochastic model (Definition~\ref{def:stoch})
where, as in the basic
 model~\eqref{eq:basicmodel}
there is an interplay between free virus, infected hepatocytes and noninfected hepatocytes.
The interactions are random and so are replications within the infected cells: with a small probability,
whenever an infected cell produces a new virion, there might be a mutation, that is a random error
in the replicative process. In particular
when a cell is infected with the type 1 virus, it is possible that it also produces type 2 virions 
(more precisely we suppose that a mutation from type 1 leads to type 2 virions),
and, still with a small probability, type 3 is generated by mutation from hepatocytes infected with the
type 2 variant. Thanks to the fact that the number of hepatocytes is very large, we are able to
approximate our stochastic model with a deterministic model which, as the basic model, is the solution
of an ODEs system, the system in \eqref{eq:diffsyst}. In Section~\ref{sec:stochmodel} we recall the 
mathematical results which allow us to approximate the average behaviour of our stochastic model
with this deterministic model (Definition~\ref{def:density} and Theorem~\ref{th:detapprox}) 
and explain how to apply them to our specific case.

In Section~\ref{sec:andet} we analyze the equilibrium points of the ODEs system \eqref{eq:diffsyst}
(the system has no exact solutions). There are four equilibria: the disease-free equilibrium $X_0$,
an equilibrium $X_1$ with all the three types of infection, $X_2$ with type 2 and type 3 infection
 together and
a fourth equilibrium $X_3$ with only the type 3 infection.
The stability of these equilibria depends on the reproductive ratios of the three types
(defined in \eqref{eq:rratios}). In particular we are able to say that a variant tends to disappear
if its reproductive ratio is smaller than 1 or if it is smaller than the reproductive ratio
of another variant which has reproductive ratio larger than 1.
We discuss how drug-resistance is explained by our model.
More precisely, therapy lowers the reproductive ratios of the affected variants and thus,
under therapy, the system moves from a situation $X_1$
 to other equilibria where the cured variants disappear and others emerge.
The main assumption is that, even if at $X_1$ all three variants are present, nevertheless type
2 and type 3 are undetectable or at least they represent only a small fraction of the total viral load,
hence they are not observed until the system reaches $X_2$ or $X_3$ (which are the ``drug-resistant states'').

In Section~\ref{sec:numerical} we focus on the main question of the paper, that is whether drug resistance
appears first under combination therapy or under sequential monotherapy. The answer relies on the numerical
solution of the ODEs system \eqref{eq:diffsyst} where therapy is considered in \eqref{eq:F} in the
appropriate time intervals (see Remark~\ref{rem:therapy}).
To compute the  numerical solutions we pick some of the unknown parameters from the literature
and deduce the missing ones as functions of the serum viral loads at equilibrium (in $X_1$, $X_2$ and $X_3$) 
of the three variants. Letting these viral loads vary among several possible configurations
(which depend on the fitness of the unknown variants type 2 and type 3) we study 46 plausible parameter sets
and the corresponding numerical solution. Each solution is tested with combination therapy and
sequential monotherapy. The time when type 3 reaches detectability under the two therapeutic approaches 
is computed. A comparison is presented, with different initial conditions.

Section~\ref{sec:conclusions} is devoted to the discussion of the consequences
of our results onto the choice of the best therapeutic approach in the cure of chronic HBV infection.
In particular our results suggest that combination therapy is the best choice if started at the early
stages of the chronic HBV infection.
 
%
%
%



\section{Stochastic model and deterministic approximation}
\label{sec:stochmodel}

The simplest model for viral reproduction is the branching process
(\cite{cf:Harris63}). This model is very rough (the behaviour of the average number of
particles is either exponential growth or exponential decrease), therefore to add complexity,
one introduces space and/or interaction between particles, see \cite{cf:DurrLevin}.
Adding space one obtains the branching random walk where the reproductive capacity of the virus 
depends on its ``location'' (which may represent not only position but also type):
its behaviour on complex networks has been studied for instance in \cite{cf:BZ,cf:BZ2,cf:Z1}.
In order to obtain a more realistic behaviour, where only a certain maximal viral load can
be achieved (or at least, growth slows down when the viral population is too large),
we can put interaction into play (logistic competition is perhaps the simplest interaction).
Spatially displaced interacting particle systems have been studied for instance in
\cite{cf:BBZ,cf:BLZ,cf:BPZ}. 

Here we present a model, with no space, which is essentially a modification of the 
basic model \eqref{eq:basicmodel}, is stochastic and describes more complex interactions between
the involved ``particles''.
The actors in the basic model are $T$, $Y$ and $V$,
representing the total number of uninfected hepatocytes, of infected hepatocytes
and of free virions, respectively. We first justify a modification of this model and introduce
a stochastic description, then put different genomic variants into play.

We assume that, by homeostasis, the total number of hepatocytes
is constant, that is $T+Y=H_0$, where $H_0$ is fixed.
Namely we suppose that, whenever an infected hepatocyte dies, 
it is replaced (almost immediately) by an uninfected one. Of course this means that infected hepatocytes
are not destroyed by the immune system at a too fast rate.
This is in agreement with the common belief that cytotoxic T lymphocytes may 
directly inhibit viral replication and thus inactivate HBV without killing the infected hepatocyte,
and that in chronically infected patients, infected hepatocytes escape immune recognition, sparing
massive liver damage (see \cite{cf:Cytotoxic,cf:tokill,cf:HepB}).
This is particulary true during the first stage of the disease, the compensated phase (which lasts 
years under a correct antiviral therapy) or in the large group of HBV infected people called 
\textit{inactive carriers} \cite{cf:Sharma}.
This assumption
allows us to deal only with the dynamics of $Y$ and $V$ and to replace $T$ with $H_0-Y$.
Since the numbers $Y$ and $V$ change according to the random interactions between free virus,
hepatocytes (infected or not) and immune system, it is natural to treat 
$Y(t)$ and $V(t)$ as random variables. More precisely, 
 $(Y(t),V(t))_{t\ge0}$ is a continuous time Markov
chain, with values in $\N^2$ where $\N$ is the set of nonnegative integers (including zero). 
The random interactions are as follows
\begin{enumerate}[(i)]
 \item a single virion infects target cells at a rate which is proportional 
to the fraction of target hepatocytes among the total, namely $\alpha \frac{H_0-Y}{H_0}$;
 \item a single infected hepatocyte dies at rate $d$;
 \item virions are produced by any infected cell at rate $p$;
 \item virions are cleared from bloodstream at rate $c$.
\end{enumerate}
\begin{figure}[h]
\begin{tikzpicture}[scale=2]
 \draw (0, 0) ellipse (.85 and 0.5);
 \node at (0,0) {$T$};
 \draw (3, 0) ellipse (.85 and 0.5);
 \node at (3,0) {$Y$};
 \draw[->] (2.5,0.5) .. controls (1.5,0.7) .. (0.5, 0.5);
\node at (1.5,0.8) {$\delta$};
\draw[->] (1,0)--(2,0);
\node at (1.5,0.1) {$\alpha$};
\draw[->] (2.15,-0.1)--(1.9,-0.3);
\node at (2,-0.45) {$p$};

\node at (1.1,-0.4) {$V$};
\draw (1.3, -0.3) circle (.03);
\draw (1.5, -0.3) circle (.03);
\draw (1.7, -0.3) circle (.03);

\draw (1.4, -0.2) circle (.03);
\draw (1.6, -0.2) circle (.03);

\draw (1.4, -0.4) circle (.03);
\draw (1.6, -0.4) circle (.03);

\draw[->] (1.5,-0.5)--(1.5,-0.75);
\node at (1.6,-0.65) {$c$};
\node at (1.4,-0.85) {out of the system};
\end{tikzpicture}
\caption{The dynamics of the one-type system.}
\end{figure}
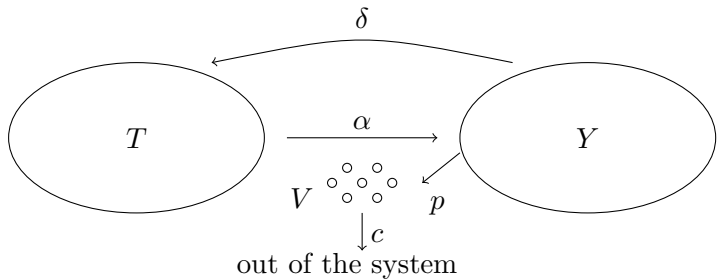
Thus the dynamics of $(Y(t),V(t))_{t\ge0}$ is described by the following rates 
($q({ k, l})$ 
denotes the transition rate from $k$ to $l$, $k,l\in\N^2$):
\[\begin{cases}
   q({(Y,V),(Y+1,V)}) & = \alpha \frac{H_0-Y}{H_0}V;\\
   q({(Y,V),(Y-1,V)}) & =  \delta Y, \text{ if }Y\ge1;\\
   q({(Y,V),(Y,V+1)}) & = pY;\\
    q({(Y,V),(Y,V-1)}) & = cV, \text{ if }V\ge1;\\
    q({ k, l}) & = 0\text{ otherwise.}
  \end{cases}
\]

Now we introduce the three different variants of the HBV: the pre-existing variant type 1;
type 2 which appears after therapy with drug A as a mutant from type 1;
type 3 which is a mutation from type 2, selected by therapy with drug B.
We denote by $V_1$ the total number of type 1 virions, $V_2$ the total number of type 2
virions, $V_3$ the total number of type 3 virions, and by $Y_1$, $Y_2$, $Y_3$
the number of hepatocytes that are carrying the corresponding type of infection
(we suppose that an hepatocyte can be infected only by one virion at a time).
The process is characterized by the following positive parameters: 
\begin{enumerate}[(a)]
 \item the infection parameters
$\alpha_1$, $\alpha_2$, $\alpha_3$, where $\alpha_i$ corresponds to the rate at which virions of
type $i$ infect target hepatocytes ($i=1,2,3$);
\item
 the production parameters $p_1$, $q_2$, $r_3$, i.e.~the rates at which hepatocytes which are 
infected with type 1, type 2 and type 3 virus respectively,
produce virions of the same type;
\item
the mutation parameters $p_2$ and $q_3$,  i.e.~the rates at which hepatocytes 
which are infected with type 1 and type 2 virus respectively, produce virions of type 2 and type 3;
\item
the hepatocytes clearance rates $\delta_1$, $\delta_2$, $\delta_3$, where $\delta_i$ corresponds 
to the rate at which hepatocytes carrying type $i$ infection are removed and thus replaced by target 
hepatocytes ($i=1,2,3$);
\item
the free virions clearance rates $c_1$, $c_2$, $c_3$, where $c_i$ corresponds 
to the rate at which virions of type $i$ are cleared from bloodstream ($i=1,2,3$).
\end{enumerate}

\begin{defn}[Stochastic model]\label{def:stoch}
Given the positive parameters $\alpha_i$, $\delta_i$, $c_i$ ($i=1,2,3$) and $p_1$, $q_2$, $r_3$, $p_2$, $q_3$,
the stochastic process  $(Y_1(t),Y_2(t),Y_3(t),V_1(t),V_2(t),V_3(t))_{t\ge0}$ is a continuous time
 Markov chain with values in $\N^6$. The transition rates $q_{ k, k+l}$ ($k\in\N^6$, $l\in\Z^6$) are as follows, 
for all $ k\neq\mathbf0$,
 \[\begin{cases}
   \alpha_i \left(1-\frac{Y_1}{H_0}-\frac{Y_2}{H_0}-\frac{Y_3}{H_0}\right)
	V_i& \text{ if } l=\mathbf e_i, i=1,2,3;\\
   \delta_iY_i &\text{ if }l=-\mathbf e_i
	      ,i=1,2,3;\\
    p_1Y_1&\text{ if }l=\mathbf e_4;\\
    p_2Y_1+q_2Y_2&\text{ if }l=\mathbf e_5;\\
    q_3Y_2+r_3Y_3&\text{ if }l=\mathbf e_6;\\
    c_iV_i & \text{ if }
     i=1,2,3\text{ and }l=-e_4,-e_5,-e_6,\text{ respectively};\\
0 & \text{ otherwise};
  \end{cases}
\]
($\mathbf e_i$ being the elements of the natural base of $\R^6$).
\end{defn}
Although this is a random process, results by Kurtz prove that its dynamic behaviour is
approximated by the solution of an ODEs system
(\cite{cf:Kurtz1}, \cite{cf:Kurtz2},
see also \cite{cf:Sani} for an example of application in 
epydemics models).
These results apply to density dependent processes, hence we must consider a modification of
our process.
\begin{defn}[Density dependent process and density process]\label{def:density}
A one-parameter family of continuous time Markov chains $(X^{(N)}(t))_{t\ge0}$
with state space $E\subseteq\Zd$ and
transition rates $(q_{ij})$ is called density dependent if there exists a continuous function
$f :\Rd \times\Zd\to \R$, such that
\[
 q_{k,k+l}=Nf\left(\frac kN,l\right),\qquad l\neq0 \text{ and } k,l\in\Zd.
\]
Suppose $(X^{(N)}(t))_{t\ge0}$ is a density dependent process.
By rescaling with $N$ we get the density process 
$(X_N(t))_{t \ge0}:=(\frac1NX^{(N)}(t))_{t \ge0}$.
\end{defn}
Under certain conditions $(X_N(t))_{t\ge0}$ converges to a deterministic process that is
the solution of a system of first order ODEs that is governed by the following function $F$:
\[
 F(x)=\sum_{l\in\Zd}lf(x,l),
\]
as stated in \cite[Theorem 3.1]{cf:Kurtz1}. We recall the theorem here.
\begin{teo}[Deterministic Approximation]\label{th:detapprox}
  Suppose that there exists an open set $E \subseteq \Rd$
the function $F$ is Lipschitz 
continuous and
\[
 \sup_{x\in E}\sum_l |l|f(x,l)<\infty;\qquad 
      \lim_{r\to\infty}\sup_{x\in E}\sum_{|l|>r} |l|f(x,l)=0.
\]
Then, for every trajectory $(x(s,x_0), s \ge0)$ satisfying the following system of ODEs
\[
 \begin{split}
  \frac d{ds}x(s,x_0)&=F(x(s,x_0));\\
  x(0,x_0)&=0, \quad x(s,x_0)\in E, 0\le s\le t,
 \end{split}
\]
$\lim_{N\to\infty}X_N(0)=x_0$ implies that for every $\delta>0$,
\[
\lim_{ N\to\infty}\Prob(\sup_{0\le s\le t}|X_N(s)-x(s,x_0)|>\delta)=0.
\]
\end{teo}
This theorem implies that the process $(X_N(t))_{t\ge0}$ can be approximated to first order by a deterministic
process, for large N. If the density process $(X_N(t))_{t\ge0}$ is initially close to $x_0$, it will tend to stay
close to the trajectory $(x(s,x_0), s > t)$ in some appropriate time-interval, subject to small random
oscillations about the path.
It is even possible to describe the behaviour of the random fluctuations of the density process
$(X_N(t))_{t\ge0}$ around its deterministic approximation (\cite[Ch.11]{cf:Kurtz}, see also \cite{cf:Sani}),
but this goes beyond the aim of this paper.

Now we apply these results to the stochastic model defined in Definition~\ref{def:stoch}.
If we take $N=H_0$, we have
 \[f( x, l)=\begin{cases}
   \alpha_i \left(1-x_1-x_2-x_3\right)x_{i+3}& \text{ if } l=\mathbf e_i, i=1,2,3;\\
   \delta_ix_i &\text{ if }l=-\mathbf e_i
	      ,i=1,2,3;\\
    p_1x_1&\text{ if }l=\mathbf e_4;\\
    p_2x_1+q_2x_2&\text{ if }l=\mathbf e_5;\\
    q_3x_2+r_3x_3&\text{ if }l=\mathbf e_6;\\
    cx_i & \text{ if }l=-\mathbf e_i, i=4,5,6;\\
0 & \text{ otherwise}.
  \end{cases}
\]
Thus $F$ has the following expression and it is easily checked that the 
hypotheses of Theorem~\ref{th:detapprox} are satisfied:
\begin{equation}\label{eq:F}
  F(y_1,y_2,y_3,v_1,v_2,v_3)=\begin{pmatrix}
       \alpha_1 \left(1-y_1-y_2-y_3\right)v_{1}-\delta_1y_1\\
       \alpha_2 \left(1-y_1-y_2-y_3\right)v_{2}-\delta_2y_2\\
       \alpha_3 \left(1-y_1-y_2-y_3\right)v_{3}-\delta_2y_3\\
       p_1y_1-cv_1\\
       p_2y_1+q_2y_2-cv_2\\
       q_3y_2+r_3y_3-cv_3
      \end{pmatrix}.
\end{equation}
Since $H_0=2\cdot 10^{11}$
(\cite{cf:liverbook}, also in accordance with the result of 139$\cdot10^6$ cells/g 
in \cite{cf:hepatocytes} and an average weight for a human liver of 1.5$\cdot10^3$g),
we may infer that the trajectories of the solution of the following system
of ODEs
\begin{equation}\label{eq:diffsyst}
\left(
 \dot y_1, \dot y_2, \dot y_3, \dot v_1, \dot v_2, \dot v_3
\right)^T
= F(y_1,y_2,y_3,v_1,v_2,v_3).
\end{equation}
are a good approximation of the ones of the stochastic process 
$\left(\frac{Y_1(t)}{H_0},\frac{Y_2(t)}{H_0},\frac{Y_3(t)}{H_0},
\frac{V_1(t)}{H_0},\frac{V_2(t)}{H_0},\frac{V_3(t)}{H_0}
\right)_{t\ge0}$,
at least in a time interval $[0,t_0]$ (where we assume that $t_0$ is sufficiently large
for our purposes).
This means that $y_1(t)$, $y_2(t)$, $y_3(t)$ denote the fraction of hepatocytes that
are infected, at time $t$, with type 1, type 2 and type 3 virus respectively.
The total number of virions of type $i$ ($i=1,2,3$), flowing in the bloodstream at time $t$, is
approximated by $v_i(t)\cdot 2\cdot 10^{11}$ and the corresponding serum viral load is obtained
dividing this quantity by $3\cdot 10^3$ (the average number of ml of serum in the human body).
\begin{rem}\label{rem:therapy}
The stochastic model we defined in \eqref{def:stoch} and its approximation, that is the solution
of \eqref{eq:diffsyst} are not only models for the case where there is no therapy, but also 
for the ones where there is therapy with drug A or B or both. We only have to change the parameters
involved by therapy. To be precise, in the first case we must 
replace $\alpha_1$, $p_1$ and $p_2$ by
 $(1-\eta_A)\alpha_1$, $(1-\varepsilon_A)p_1$ and $(1-\varepsilon_A)p_2$ respectively; 
in the second and third case $\alpha_2$, $q_2$ and $q_3$ by  $(1-\eta_B)\alpha_2$, $(1-\varepsilon_B)q_2$ 
and $(1-\varepsilon_B)q_3$ respectively; in the third case $\alpha_1$, $p_1$ and $p_2$ must be replaced by
 $(1-\eta_C)\alpha_1$, $(1-\varepsilon_C)p_1$ and $(1-\varepsilon_C)p_2$ respectively, where $\eta_C$
and $\varepsilon_C$ are the efficacy parameters of the combination of drug A and B against the type 1 
infection.
\end{rem}

\section{Analysis of the deterministic model}\label{sec:andet}

The ODEs system \eqref{eq:diffsyst} has no exact solution, thus in order to understand its behaviour,
we start by analyzing its equilibria.
Solving $F(x)=0$ ($F$ taken from \eqref{eq:F}), we find that there are four equilibrium points, 
the disease-free state and three infected equilibria. 
The Jacobian matrix of \eqref{eq:F}
is $J(x)$ given by 
\[
\begin{pmatrix}
       -(\alpha_1v_1+\delta_1) & -\alpha_1v_1 &  -\alpha_1v_1 & \alpha_1\beta & 0 & 0\\
       -\alpha_2v_2 &  -(\alpha_2v_2+\delta_2) &  -\alpha_2v_2 & 0& \alpha_2\beta  & 0\\
        -\alpha_3v_3 &  -\alpha_3v_3 & -(\alpha_3v_3+\delta_3) &  0& 0 & \alpha_3\beta \\
         p_1 & 0 & 0 & -c_1 & 0 & 0 \\
         p_2 & q_2 & 0 & 0 & -c_2 & 0\\
         0 & q_3 & r_3 & 0 & 0 & -c_3
      \end{pmatrix},
\]
where $\beta=1-y_1-y_2-y_3$.

Before going into the analysis of the equilibria,
let us mention a few words about the choice of the parameters and what we expect to
observe. The viral clearance parameters $c_i$ depend only on the life-time of virions in bloodstream,
thus on the antibodies production of the infected patient. We assume that this production does
 not depend on the variant, thus from now on we put $c:=c_1=c_2=c_3$.
We also assume that $p_2\ll p_1$ and $q_3\ll q_2$. Indeed $p_2$ and $q_3$ represent the
reproduction rate of type 1 and type 2 infected hepatocytes respectively, multiplied by the probability
that a mutation appears (from type 1 to type 2 and from type 2 to type 3 respectively).
The probability of mutation is very small, thus we assume that not only
$p_2/p_1$ and $q_3/q_2$ are very small but also $q_3$ and $p_2/q_2$ are small
(in Section~\ref{sec:numerical} we will choose $q_3=10^{-7}$ or $q_3=10^{-6}$
and computation will lead us to, in different cases, $p_2/q_2$ between $10^{-10}$ and $10^{-4}$).
Therefore in the analysis of some of the equilibria, we will approximate $p_2/q_2$ and $q_3$ with zero.
An important role in the stability of the equilibria is played by the basic reproductive
ratios of the three variants
\begin{equation}\label{eq:rratios}
  R_1=\frac{\alpha_1p_1}{c\delta_1},\qquad 
 R_2=\frac{\alpha_2q_2}{c\delta_2},\qquad
 R_3=\frac{\alpha_3r_3}{c\delta_3}.
\end{equation}
Note that if the patient is undergoing therapy the corresponding infection and production
parameters have to be changed according to Remark~\ref{rem:therapy}.
For instance, under therapy with drug A, the basic reproductive ratio of the type 1 variant becomes
$(\alpha_1p_1(1-\eta_A)(1-\varepsilon_A))/(c\delta_1)$.

\subsection{Disease-free equilibrium}
The disease-free equilibrium is $X_0=(0,0,0,0,0,0)$. 
Evaluating at $X_0$ we get
\[
 J(X_0)=\begin{pmatrix}
       -\delta_1 & 0 &  0 & \alpha_1 & 0 & 0\\
        0 &-\delta_2 &   0 & 0& \alpha_2  & 0\\
        0 &  0 & -\delta_3 &  0& 0 & \alpha_3 \\
         p_1 & 0 & 0 & -c & 0 & 0 \\
         p_2 & q_2 & 0 & 0 & -c & 0\\
         0 & q_3 & r_3 & 0 & 0 & -c
      \end{pmatrix}.
\]
This matrix has six eigenvalues:
\[
  \frac{-\delta_1-c\pm\sqrt{4\alpha_1p_1+(\delta_1-c)^2}}{2},
 \frac{-\delta_2-c\pm\sqrt{4\alpha_2q_2+(\delta_2-c)^2}}2 ,
\frac{-\delta_3-c\pm\sqrt{4\alpha_3r_3+(\delta_3-c)^2}}{2}.
\]
These eigenvalues are real and they are all negative (hence $X_0$ is locally asymptotically stable)
if $R_i<1$ for all $i=1,2,3$. If at least one $R_i$ is larger than 1, then there is at least
one eigenvalue which is positive and $X_0$ is locally asymptotically unstable.

\subsection{The infected equilibrium without type 1 and type 2 infections.}
The infected equilibrium without type 1 infection is $X_3=(0,0,\tilde y_3,0,0,\tilde v_3)$
and has biological meaning (that is, $\tilde y_3$ and $\tilde v_3$ are positive) only if $R_3>1$, since
\[
 \tilde v_3=(R_3-1)\delta_3/\alpha_3,\qquad
 \tilde y_3=1-1/R_3=\alpha_3\tilde v_3/\delta_3R_3.
\]
The Jacobian matrix, evaluated at $X_1$, has six eigenvalues
\[\begin{split}
 &\frac{-\alpha_3r_3(c+\delta_1)\pm\sqrt{(\alpha_3r_3(c-\delta_1))^2+4\alpha_3r_3\alpha_1p_1c\delta_3
      }}{2\alpha_3r_3} ,\\
&\frac{-\alpha_3r_3(c+\delta_2)\pm\sqrt{(\alpha_3r_3(c-\delta_2))^2+4\alpha_3r_3\alpha_2q_2c\delta_3
      }}{2\alpha_3r_3} ,\\
 &\frac{-\alpha_3r_3-c^2\pm\sqrt{4c^3\delta_3+(\alpha_3r_3-c^2)^2}}{2c} .
  \end{split}
\]
They are all real and easy computations show that $X_3$ is locally asymptotically 
stable if $R_3>1$, $R_3>R_1$ and  $R_3>R_2$. If $R_3<1$ or
$R_2>R_3$ or  $R_1>R_3$, then $X_3$ is locally asymptotically 
unstable.

\subsection{The infected equilibrium without type 1 infection.}
The infected equilibrium without type 1 infection is $X_2=(0,\bar y_2,\bar y_3,0,\bar v_2,\bar v_3)$
and has biological meaning only if $R_2>1$ and $R_2>R_3$, since
\[\begin{array}{ll}
   \bar v_2= \frac{q_2}{c}\frac{(R_2-1)(R_2-R_3)}{R_2(R_2-R_3+R_3q_3/r_3)},\qquad
&\bar v_3=\frac{q_3(R_2-1)}{c(R_2-R_3+R_3q_3/r_3)}=\frac{q_3}{q_2}\frac{R_2}{R_2-R_3}\bar v_2\\
\bar y_2=\frac{c}{q_2}\bar v_2, 
&\bar y_3=\frac{c}{r_3}\frac{R_3}{R_2}\bar v_3.
  \end{array}
\]
Two eigenvalues are
\[
 \frac12\left(-c-\delta_1\pm {\sqrt{(\delta_1-c)^2+4\alpha_1p_1(1-\bar y_2-\bar y_3)}}\right),
\]
which are negative if and only if $R_1<1/(1-\bar y_2-\bar y_3)= R_2$.
The other four eigenvalues are the solutions of a quartic equation and have expressions
which are too involved to be studied here. We prefer to study the approximated
solutions obtained substituting $q_3=0$ (and $\bar v_3\approx0$, which is plausible if $R_2-R_3$
is not too small) in the characteristic polynomial:
\[\begin{split}
   &\frac12\left(-\alpha_2v_2-\delta_2-c\pm
   \sqrt{(\alpha_2\bar v_2+\delta_2-c)^2+4\alpha_2q_2(1-\bar y_2-\bar y_3)}\right)\\
    &\frac12\left(-c-\delta_3\pm \sqrt{(\delta_3-c)^2+4\alpha_3r_3(1-\bar y_2-\bar y_3)}
    \right).
  \end{split}
\]
Computations show that these approximate eigenvalues are all negative if $R_2>R_3$
and $\bar v_2>0$. Thus $X_2$ is locally asymptotically stable if $R_2>1$ and $R_2>R_3$. 
If at least one of these conditions
is not satisfied, $X_2$ is locally asymptotically unstable.

\subsection{The infected equilibrium with type 1}
The infected equilibrium where also type 1 is present,
is $X_1=(\hat y_1,\hat y_2,\hat y_3,\hat v_1,\hat v_2,\hat v_3)$
where 
\begin{equation}\label{eq:X3}
 \begin{array}{ll}
  \hat v_1=\frac{p_1}{c}\hat y_1, \qquad &\hat v_2=\frac{R_1}{R_2}\hat y_2,
\qquad\hat v_3= \frac{R_1}{R_3}\frac{r_3}{c}\hat y_3,\\
\hat y_2=\frac{R_2(p_2/q_2)}{R_1-R_2}
  \hat y_1, & \hat y_3=\frac{p_2q_3}{q_2r_3}\frac{R_2R_3}{(R_1-R_2)(R_1-R_3)}.
 \end{array}
\end{equation}
As for
$\hat y_1$, it has a complicated expression which is of no interest
here (but can be obtained from the relation $\hat y_1+\hat y_2+\hat y_3=1-1/R_1$). Here is its
approximation if one substitutes $q_3\approx0$:
\[
    \hat y_1\approx 
 \frac{(R_1-1)(R_1-R_2)}{R_1(R_1-R_2+R_2(p_2/q_2))},
\]
thus $X_1$ has biological meaning if $R_1>1$, $R_1>R_2$ and $R_1>R_3$.
The eigenvalues of $J(X_1)$ have too complicated expressions, hence we prefer to study the approximate
solutions obtained with the substitutions $q_3=0$, $p_2=0$ (and $\hat v_2\approx0$, 
which is plausible if $p_2/q_2$ is small and $R_1-R_2$
is not too small) in the characteristic polynomial.
The approximated eigenvalues are
\[\begin{split}
   &\frac12\left(-\alpha_1\hat v_1-\delta_1-c\pm
   \sqrt{(\alpha_1\hat v_1+\delta_1-c)^2+4\alpha_1p_1(1-\hat y_1-\hat y_2-\hat y_3)}\right)\\
   &\frac12\left(-c-\delta_2\pm
   \sqrt{(\delta_2-c)^2+4\alpha_2q_2(1-\hat y_1-\hat y_2-\hat y_3)}\right)\\
    &\frac12\left(-c-\delta_3\pm \sqrt{(\delta_3-c)^2+4\alpha_3r_3(1-\hat y_1-\hat y_2-\hat y_3)}
    \right).
  \end{split}
\]
Substituting $\hat y_1+\hat y_2+\hat y_3= 1-1/R_1$,
we get that $X_1$ is locally asymptotically stable if $R_1>1$, $R_1>R_2$ and $R_1>R_3$
and unstable otherwise.

The conditions for stability can be summarized by Table 1. 
\begin{table}[h]\label{tab:table1}
 \begin{tabular}{|l|l|l|l|}
 \hline
$R_1$ & $R_2$ & $R_3$ & \textbf{Equilibria}\\
\hline
$R_1<1$ & $R_2<1$ & $R_3<1$ & $X_0$ stable, other $X_i$s unstable\\
\hline
$R_1>1$ & $R_2<R_1$ & $R_3<R_1$ & $X_1$ stable, other $X_i$s unstable\\
\hline
$R_1<R_2$ & $R_2>1$ & $R_3<R_2$ & $X_2$ stable, other $X_i$s unstable\\
\hline
$R_1<R_3$ & $R_2<R_3$ & $R_3>1$ & $X_3$ stable, other $X_i$s unstable\\
\hline
\end{tabular}
\bigskip
\caption{Stability of the equilibria as a function of $R_1$, $R_2$, $R_3$.}
\end{table}
\subsection{Drug resistance explained by the model.}
Let us discuss how the presence of these equilibria and their stability can be biologically interpreted.
We are modelling the infection within a chronic patient, therefore with no treatment all three
variants have reproductive ratios larger than 1. Moreover we believe that type 1 has largest fitness
and type 2 has larger fitness than type 3. Therefore in the beginning $R_1>R_2>R_3>1$.
With no cure, after an appropriate amount of time, the system tends to lie in a neighbourhood of $X_1$.
We believe that the reason why drug resistant mutants are not observed unless under therapy
is the fact that, even if mutations do happen, nevertheless in the equilibrium together with
the wild type, mutants are numerically negligible.
In other words, if we denote by $v_{i,j}$ the $v_i$-th component of the equilibrium point $X_j$,
we assume that $v_{2,1}$ and $v_{3,1}$ are negligible if compared to $v_{1,1}$ (we already used
this assumption in the analysis of equilibria).
If the patient is given drug A, then type 1 is cured ($R_1$ is lowered, by multiplication
by $(1-\eta_A)(1-\varepsilon_A)$), then the system moves towards equilibrium $X_2$, where
variant 2 emerges.  If also type 2 is cured, with drug B ($R_2$ is lowered, by multiplication
by $(1-\eta_B)(1-\varepsilon_B)$), the new equilibrium will be $X_3$, that is, 
the patient will show drug resistance. 

Studies show that combination therapy has a larger efficacy in reducing the production rates
(compared to single drug therapy); 
for instance \cite{cf:lau-tsiang} showed that, at the dosage of their trial,
lamivudine only has $\varepsilon=0.94$, while combined with famciclovir it gets $\varepsilon=0.988$.
We suppose that whenever a virus is the target of 
two drugs together, the effects add in
an independent manner, namely the combined efficacies of the cocktail drug A + drug B
against production of type 1 and against infection of new cells respectively, are 
$\varepsilon_C=1-(1-\varepsilon_A)(1-\varepsilon_B)=
\varepsilon_A+\varepsilon_B-\varepsilon_A\varepsilon_B$ and
$\eta_C=1-(1-\eta_A)(1-\eta_B)=
\eta_A+\eta_B-\eta_A\eta_B$. 
If the patient is given drug A and drug B in combination, then the viral load of type 1
decays more rapidly ($R_1$ is multiplied by $(1-\eta_C)(1-\varepsilon_C)$)
and the system moves directly from $X_1$ to $X_3$.
The disease-free equilibrium represent the desirable state where the patient is permanently cured,
and can be reached only if we have a drug at hand which acts also against the third variant
(by assumptions, this is not the case under study).
We believe that the fact that type 1 has the 
largest fitness is also reflected in the fact
that when type 1 is removed from the system, the other two types can never exceed type 1's initial
viral load. Moreover, viral loads of a fixed type increase when the competitors are erased.
Based on these observations, we make some more assumptions on the $v_{i,j}$s: 
$v_{1,1}\gg v_{2,1}\ge v_{3,1}$, $v_{1,1}\ge v_{2,2}\ge
v_{3,3}> v_{3,2}\ge v_{3,1}$,
$v_{2,2} > v_{2,1}$.

The qualitative behaviour in time of the logarithm of viral loads under combination or sequential therapy 
(if therapy starts when the system is at $X_1$) can be described
by Figure 2 (where $v_1(t)$ is represented by a black line, $v_2(t)$ by a dashed line and
$v_3(t)$ by the dotted one). Note that the picture is simplified, it does not show the biphasic decrease of
viral loads under therapy, which is observed in trials and has been explained in \cite{cf:colombatto}.
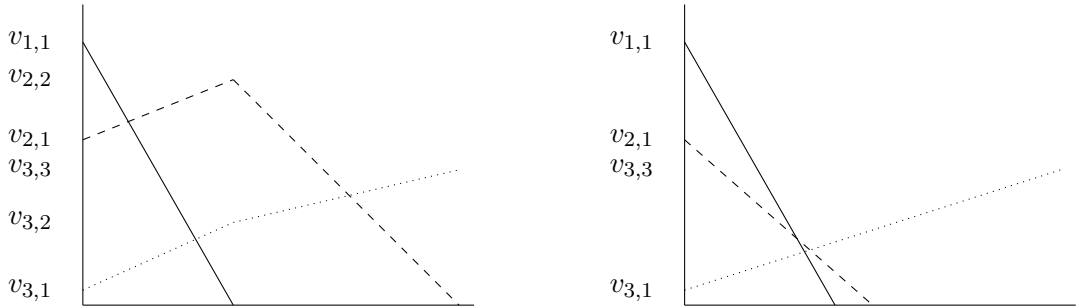
\begin{figure}[ht]
 \begin{tikzpicture}
 \draw (0,0)--(5.2,0);
 \draw (0,0)--(0,4);
 \node  at (-.7,3.5) {$v_{1,1}$};
 \node  at (-.7,2.2) {$v_{2,1}$};
 \node  at (-.7,.2) {$v_{3,1}$};
 \node  at (-.7,3) {$v_{2,2}$};
 \node  at (-.7,1.1) {$v_{3,2}$};
 \node  at (-.7,1.8) {$v_{3,3}$};

 \draw (0,3.5)--(2,0);
 \draw[style=dashed] (0,2.2)--(2,3);
 \draw[style=dashed] (2,3)--(5,0);
 \draw[style=dotted] (0,.2)--(2,1.1);
 \draw[style=dotted] (2,1.1)--(5,1.8);
 \draw (8,0)--(13.2,0);
 \draw (8,0)--(8,4);
 \node  at (7.3,3.5) {$v_{1,1}$};
 \node  at (7.3,2.2) {$v_{2,1}$};
 \node  at (7.3,.2) {$v_{3,1}$};
 \node  at (7.3,1.8) {$v_{3,3}$};

 \draw (8,3.5)--(10,0);
 \draw[style=dashed] (8,2.2)--(10.5,0);
 \draw[style=dotted] (8,.2)--(13,1.8);

\end{tikzpicture}
\caption{Viral loads dynamics under sequential (left) and combination (right) therapy.}
\end{figure}

\section{Numerical comparison under sequential or combination therapy}
\label{sec:numerical}


The best therapeutic approach will be the one which sets off for a longer time the appearance
of a detectable viral load of type 3 virions in the bloodstream.
It is clear that the very presence of type 1 and type 2 virions in the bloodstream and in infected
hepatocytes is an antagonist to the development of type 3, since all the variants are competing
for a finite resource, namely, the target, non-infected hepatocytes. On the other hand, whenever
type 1 (type 2 respectively) replicates, there is a chance of producing new variants (type 2 and type 3
respectively) that may be capable of drug resistance. Thus, in principle, it is not trivial to
answer to the question whether it is a better strategy to employ sequential monotherapy
(which may lead to surface type 2 virions and afterwards type 3 virions)
 or combination therapy (which strikes down both type
1 and type 2 but may clear the way to type 3).
The answer relies on the quantitative dependence on time of viral loads.
In order to analyze it, we resort to numerical solutions of the ODEs system \eqref{eq:diffsyst}.

We need to know the numerical values
of the parameters involved, and of the initial condition. 
The twelve parameters are:
the infection parameters $\alpha_1$, $\alpha_2$ and $\alpha_3$; the clearance 
parameters $c$, $\delta_1$, $\delta_2$ and $\delta_3$;
the production parameters $p_1$, $p_2$, $q_2$, $q_3$, $r_3$.
The initial condition is given by  $(y_1(0), y_2(0) , y_3(0) ,v_1(0) ,v_2(0) , v_3(0))$. 
Moreover, we also need the efficacy parameters $\eta$ and $\varepsilon$
of the drugs involved in the testing. We gather from the literature 
estimates on some of the parameters regarding the wild-type chronic infection and on the efficacy
parameters.
From \cite{cf:TsiangRooney} we have 
$c=0.65$/day, $\delta_1=0.0143$/day, $\varepsilon(\text{adefovir})=0.993$, while from \cite{cf:Nowak96}
we have  $\varepsilon(\text{lamivudine})$ ranging from 0.87 to 0.96 depending on the dosage.
In their estimates, both papers use the assumption that a percentage varying between 5\% and 40\% of 
the total number of hepatocytes is productively infected in chronic HBV patients 
(the source of this assumption is \cite{cf:BianchiGudat}). Thus we take $y_1(0)=0.25$, $0.2$ and $0.3$
in different simulations.
The order of magnitude of the initial viral load in the patients of \cite{cf:TsiangRooney} ranged from 
$10^7$ to $10^9$ copies/ml: for simplicity's sake, we take $v_1(0)=1$ which corresponds to
$2\cdot 10^{11}$ free virions and a viral load of $0.67\cdot 10^{8}$ (assuming that patients
have $3\cdot10^3$ml of serum).

In our numerical tests, we 
choose $\varepsilon_A=0.97$ and $\varepsilon_B=0.85$: we assumed that the first drug is highly
efficient against the wild type, while the second drug, which acts also on type 2, and is therefore
more flexible, is less efficient (these efficacy parameters are nevertheless both in the range 
of the efficacies of known drugs).
Although some papers assume that $\eta$ is neither 0 or 1 (see for instance 
\cite{cf:Lewin,cf:ADEFOVIR-LAMIVUDINA}), there
are no explicit estimates of this parameter for the various drugs, hence in different simulations
we pick $\eta_A$ and $\eta_B$ in $\{0.25,0.5\}$
(\cite{cf:Lewin} used as a plausible value 0.5). 
Finally, it is not possible to gather information on the values
of the infection and production parameters of the variants, since ideally we would like to compare
sequential and combination therapy when facing any possible, unknown variants of the HBV virus. 
We assume that the clearance parameter of the infected hepatocytes does not depend on the type
of virus which infected the cell, that is $\delta_1=\delta_2=\delta_3=0.0143$/day.

In order to compute the missing parameters, we assume that we are given the numerical values
of the $v_i$-components of the equilibria $X_1$ 
(namely  $v_{1,1}$, $v_{2,1}$ and $v_{3,1}$), 
$X_2$ ($v_{2,2}$ and $v_{3,2}$), and $X_3$ ($v_{3,3}$). Since these six equilibrium viral loads
depend on the six missing parameters (see the expressions in in Section~\ref{sec:andet}),
 we take the inverse function and derive $\alpha_1$, $\alpha_2$, $\alpha_3$, $p_2$, $q_2$ and $r_3$.
The parameter $p_1$ is derived from $p_1=\frac{c\cdot v_{1,1}}{y_{1,1}}$, 
once we fix the value of $y_{1,1}$ ($y_{i,j}$ being the $y_i$-component of $X_j$).
We fix $y_{1,1}=y_1(0)$. We are left with an arbitrary choice of $q_3$, so we assume that it is
very small and choose it equal to $10^{-7}$.

We proceed computing all the parameter sets which are obtained when the equilibrium viral loads
vary among plausible values (we test for each value different powers of 10). 
We do not consider values that lead to nonpositive parameters
(among the conditions there is $v_{3,2}>v_{2,2}\frac{v_{3,1}}{v_{2,1}}$);
moreover we require that 
$p_1\ge q_2\ge r_3$ (type 1 has a faster production rate, and type 2 is faster than type 3).
The resulting ranges, which reflect our hypotheses on the equilibria
(for instance the negligibility of $v_{2,1}$ and $v_{3,1}$), are listed in Table 2.
\begin{table}[h]\label{tab:table2}
 \begin{tabular}{l|l|l|l|l|l}
& $v_{2,1}$  & $v_{3,1}$    & $v_{2,2}$  & $v_{3,2}$   & $v_{3,3}$ \\
\hline
$\min$ & $10^{-9}$ & $10^{-9}$ & $10^{-6}$ & $10 v_{2,2}\frac{v_{3,1}}{v_{2,1}}$ & $10^{-6}$\\
\hline
$\max$ & $\min\left(10^{-5},\frac{v_{2,2}}{10}\right)$ & $v_{2,1}$ & 1 & $\min \left(\frac{v_{2,2}}{10},\frac{v_{3,3}}{10}\right)$ & 1  \\
\hline
\end{tabular}
\bigskip
\caption{Range of equilibrium viral loads.}
\end{table}
%
We also require that 
that, at equilibria, the percentage of infected hepatocytes lies between 5\% and 40\%
(parameter sets which do not satisfy all these conditions are disregarded).
Under these conditions we get 46 cases to analyze, which correspond to different possible type 2
and type 3 variants (with more or less ability to infect and reproduce).

Once we obtain the parameters, we numerically compute the solutions of the ODEs system
in case of combination therapy
(drug 1 and drug 2 given together until type 1 and type 2 are both eradicated) or of
sequential monotherapy (drug 1 is given until type 1 is cured and we start giving drug 2 
when type 2 reaches a viral load of 667 copies/ml). 
Combination therapy is judged as the best therapeutic approach if the time when
 viral load of type 3 exceeds 20 copies/ml is significantly larger than the same time under
sequential monotherapy.

In the numerical computations the arbitrary parameters such as $\eta_A$, $\eta_B$, $y_1(0)$ and $q_3$
have been varied (one at a time) with no significant changes in our conclusions.
The only values on which the best therapeutic approach depends are the initial viral loads of type 3.
We assume that the initial viral load $v_2(0)$ is 7 or 67 copies/ml and $y_2(0)=\frac{R_2}{R_1}v_2(0)$
(that is, the fraction of initially infected hepatocytes carrying the type 2 infection is the same
function of the corresponding viral load as $y_{2,1}$ is a function of $v_{2,1}$).

We sum up the conclusions of our numerical simulations in Table 3. The answer to our main question
depends on the values of $vl_3(0)$ and $y_3(0)$, which 
are the initial serum viral load and the fraction of infected hepatocytes of type 3 respectively.
In the second column (CT best) we put the fraction of cases where the combination therapy is advantageous, and
among these cases we compute the minimum, maximum and average percentage delay 
($\min d$, $\max d$ and $\bar d$ respectively). Finally, $\bar T_s$ and 
$\bar T_c$ are the average times before drug resistance (i.e.~the detection of type 3) for 
patients under sequential or combination therapy respectively.
The table refers to data relative to the 46 cases obtained with the
following choice of the arbitrary parameters:
$q_3=10^{-7}$, $y_1(0)=0.25$, $v_2(0)=10^{-7}$, $\eta_A=\eta_B=0.5$.
We changed the value of these
parameters, one at a time, without significant changes in the conclusions, which is an indication
of robustness of our comparison.
It appears that the advantage of combination therapy is appalling if started before that type 3 is present
(first row of Table 3). If therapy starts when type 3 has a viral load of approximately 1 copy/ml but 
has not infected any hepatocyte yet (second row of Table 3), combination therapy is still the best 
choice in 38 cases among 46, but delayed drug resistance by more than 10\% of the time only in 18 cases.
Finally, if infection of type 3 has already progressed we look at the third row of Table 3.
There we chose  $y_3(0)=\frac{\alpha_3}{R_1\delta_3}v_3(0)$
(the fraction of infected hepatocytes corresponding to type 3 in equilibrium together with
the other strains, if $v_3(0)$ coincided with $\hat v_3$, see \eqref{eq:X3})
and obtained that while combination therapy is preferable in 16 cases out of 46, sequential
monotherapy is preferable in 9 cases out of 46. In the remaining 21 cases there is substantial
parity of outcomes under the two therapeutic approaches.

\begin{table}[h]\label{tab:table3}
 \begin{tabular}{l|l|l|l|l|l|l}
$vl_3(0)$, $y_3(0)$  & CT best   & $\min d$ & $\max d$  & $\bar d$   & $\bar T_s$ & $\bar T_c$\\
\hline
0, 0 & $46/46$ & 33\% & 90\% & 50\% & 6.3 years& 9.5 years\\
\hline
1, 0 & 38/46 &1\% & 56\%& 15\% & 5.1 years & 5.8 years\\
\hline
1, $>0$ &16/46 & 0.4\%& 7.7\%& 3\% & 2.87 years& 2.91 years\\
\hline
\end{tabular}
\bigskip
\caption{Comparison of therapeutic approaches depending on the initial conditions.}
\end{table}

%
%

\section{Discussion}\label{sec:conclusions}

We developed a stochastic model of HBV viral dynamics (with or without therapy) which accounts 
not only for different aspects of the interplay between the virus and the host, through some 
parameters, but also for the presence of heterogeneous strains in the viral progeny, due to
random errors in the replicative process. In particular we considered two consecutive mutations
of the primary infection.
The numerical value of the parameters in the model depend on the patient and on the viral strains.
Indeed the host is responsible of an immune response that counteracts the infection through
viral clearance from the bloodstream (humoral immunity, parameter $c$) and clearance of infected hepatocytes (cell-mediated immunity, parameters $\delta_i$). On the other hand, any virus 
has a specific capacity of infection of target cells (parameters $\alpha_i$) and of reproduction
through the infected cells (parameters $p_1$, $q_2$ and $r_3$ for the three variants, respectively).
Moreover whenever there is production of type 1 virus there is a small probability of producing
a type 2 virion and analogously, from type 2 infected hepatocytes there might be production of
type 3 virions (the mutation rates are $p_2$ and $q_3$ respectively).
We approximated the behaviour of the stochastic model with the behaviour of a deterministic model,
given by the solution of the ordinary differential equations system 
\eqref{eq:diffsyst}, which was then the object of our study.

Under the natural assumption $R_1>R_2>R_3$, based on the fact that more easily observed 
variants must have a larger fitness (reflected on the order of reproductive ratios),  
the model is able to explain drug resistance. Indeed, as we showed in Section~\ref{sec:andet},
if there is no therapy then in the long run the mutants, even if they appear by mutation,
remain numerically negligible. If the patient undergoes therapy which removes only the primary
virus, then the drug-resistant strain type 2 surfaces; while if therapy acts also against
this latter strain, in the long run the third variant, which by definition we suppose
drug-resistant to any drug at hand, appears.
Note that the model gives a qualitative description of the appearance of drug-resistant strains
even if the characteristics of these strains are unknown (the only assumption is that
mutants are negligible when together with the wild-type and in general have smaller fitness).

We used this model to compare the efficacy of combination therapy versus sequential monotherapy.
In this comparison, the best therapy is the one that mostly delays the appearance 
 of the drug resistant variant (i.e.~the time when this variant reaches a serum viral load of 20 cp/ml).
The analysis of the equilibria (Section~\ref{sec:andet}) tells us how the system behaves qualitatively
with or without therapy but does not tell us how long it takes to reach any of these equilibria.
This can be determined only solving \eqref{eq:diffsyst}, where the infection and production
parameters are multiplied
by $1-\bar\eta$ and by $1-\bar\varepsilon$ respectively, and $\bar\eta$ and $\bar\varepsilon$
depend on the therapy. For instance,
in sequential monotherapy for $t\in[0,t_1]$ therapy is drug A, for 
$t> t_1$ therapy is drug B, where $t_1$ is the time where
type 2 reaches a worrying viral load ($t_1$ depends on the solution itself).
The system  can be solved only numerically (through a mathematics software, we used Maple). 
To plug in numerical values of the parameters,
we took $c_1$, $\delta_1$ from the literature and assumed $c_1=c_2=c_3$ and 
$\delta_1=\delta_2=\delta_3$. Six unknown parameters ($\alpha_1$, $\alpha_2$, $\alpha_3$, $p_2$, 
$q_2$, $r_3$) were obtained from the expression of the equilibria, as functions of the viral loads at equilibrium of the three types. The parameter $p_1$ is a function of the equilibrium viral load
of the primary infection and of the corresponding fraction of infected hepatocytes
(which were gathered from the literature); for the last unknown parameter, $q_3$, we had to make an arbitrary choice (and put it equal to $10^{-6}$ or $10^{-7}$ in different simulations). 
We believe that, as long as $q_3$ is small, its numerical value can only rescale the 
time of emergence of the third variant (the smaller $q_3$, the longer the time) but should not 
change the kind of therapy which is most suitable.
As for the choice of the efficacy parameters of the drugs involved in the testing, we 
picked them from the literature.

We note here that under the assumption that $\delta_1=\delta_2=\delta_3$, for all plausible
choices of the equilibrium viral loads, we obtained $\alpha_1<\alpha_2<\alpha_3$. This means
that, if the mutants are not more capable of the wild-type to evade the cell-mediated immune response,
then only a greater ability of infecting target hepatocytes can lead them to detectability.
We performed some numerical calculations with some choices of $\delta_2$ and $\delta_3$
such that $\delta_1>\delta_2\ge\delta_3$ and obtained that in these cases it is possible that $\alpha_1>\alpha_2\ge\alpha_3$ (which is more intuitive, since one is lead to believe that
larger fitness also means greater infection ability).
Thus, in our model, drug-resistant mutants can reach detectibility if, compared to the wild-type,
they either have a greater infection ability or if they are better at masking themselves into
infected cells. This is plausible, since mutations may for instance affect epitopes recognized by 
cell-mediated immune response.

It is worth noting that, even if we were forced to make some arbitrary parameter choices
and the model is not fitted to any data (which is intrinsic in the problem itself since
we want to address unknown variants of the HBV), nevertheless the order of magnitude of
the time before drug-resistance that we obtained is in accordance with the time observed
in case trials. For instance in \cite[Table 1]{cf:mutations} many patients report,
with different drugs, resistance after a treatment period ranging from 3 to 6 years.

Our numerical analysis shows that combination therapy is by far the best choice if started at the early
stages of the chronic infection, while it seems that the advantage is not striking if therapy
is in act when drug-resistant mutants are already numerous and have started to significantly
infect the hepatocytes. This difference is not surprising, since if type 3 mutants are absent,
they are generated only by mutation from type 2, hence the need of lowering the number of type 2
virions and infected cells. On the other hand, if type 3 mutants are already present, in some cases
the competition with other strains might slow the proliferation of the type 3 infection.
Nevertheless,  we believe that in reality early therapy is  more convenient also for other reasons.
In our deterministic approximation of the model, there is a constant
production of mutants, in a fixed proportion of the existing replications. 
In reality those replications are random and the time $\tau$ of appearance of mutants is random.
If there are many replications $\tau$ is more likely to be small; therefore
a good strategy is also to keep replications low in order to delay the first appearance of mutants.
This randomization is not considered in the deterministic approximation of the model and gives
another indication in favour of early combination therapy.
The biological and clinical importance of this model is to show a clear advantage in favor of the 
combination strategy, as it could lead to a greater delay in the development of drug resistant 
variants, which have less fitness than the wild type (drug sensible) strain, but  are able eventually 
to escape from the therapy. From the analysis of the model we also want put emphasis on the need of 
an early antiviral therapy: since the probability of the emergence of mutation is directly 
proportional to the viral 
fitness, high levels of wild strain virus could generate more easily drug resistant variants 
(type 3 mutants). Consequently the probability of mutation, over time, is inversely proportional to 
the selective pressure exerted by the antiviral therapy.

By the way, there is also clinical evidence that the best choice is a broad spectrum early therapy: 
for instance in \cite{cf:notnaive} the authors observed that the emergence of entecavir-resistant
strains is more likely to appear in patients which were already lamivudine-resistant (rather
than in treatment-naive patients).

Concerning the type of combination therapy, ideally the best combination therapy for HBV infection
should consider the association between a \textit{high power} drug (with a laevorotatory structure) 
and a \textit{high genetic barrier} drug (with a flexible, plastic structure).
Power is the speed with which a drug causes the suppression of viral replication
hence it is proportional to the efficacy parameter $\varepsilon$ introduced in
Section~\ref{sec:numerical}). We can assume that 
a potent drug is a molecule that has a high affinity to the target protein. Such affinity, in addiction 
to an intrinsic chemical characteristic, could be an expression of a better adaptation to the 
chirality of nature of the target. Respecting the symmetry of nature, the laevogyre molecules are 
more selective and more powerful to the target, with fewer number of side effects.
The genetic barrier represents the number of mutations needed by the virus to replicate effectively 
in the presence of the drug (a higher genetic barrier indicates that there is a smaller set of mutants
that evade the action of the drug). 
 Since the action of a drug is an expression of the interaction active 
molecule -- target protein, we can assume that a drug with an high genetic barrier is a molecule that 
has a ``high flexibility'', which enables to adapt to different conformations of the target protein 
induced by punctiform genotypic mutations.
That synergism is something more than the simple sum of the effects of two drugs, and it contains 
a stronger action than monotherapy, even with entecavir or tenofovir, molecules that to date 
represent the best compromise between power and genetic barrier.
These two drugs show a low incidence of mutations
after 4 years of treatment: 1.2\% and 0\% respectively, \cite{cf:mutations} (although in principle
sooner or later drug-resistant mutants may appear with any drug, which is one argument in
favour of combination therapy even with these drugs at hand). As we 
already mentioned this dramatically changes with entecavir in non-treatment-naive patients (\cite{cf:notnaive}). Moreover recent results show the appearance of drug-resistance in
coinfected HBV-HIV patients, treated with tenofovir (\cite{cf:tenres}).
The authors identify in poor compliance to therapy the reason of this failure; it is our opinion
that this is another aspect to keep in mind. 
In conclusion, basing on the numerical evaluations carried out in the model, we suggest to 
intervene in the 
natural history of HBV infection immediately and with a broad-spectrum approach: a combination of 
powerful and high genetic barrier drugs. We believe that combination therapy is less likely to select 
resistant strains, especially in patients with a suboptimal compliance. We may infer the same 
conclusions also considering the evidence-based effectiveness of the combination therapy 
(compared to the monotherapy) in other viral infections supported by multiple subtypes, 
termed viral quasispecies, such as hepatitis C virus (HCV) and human immunodeficiency virus 
(HIV) infections.

\section*{Acknowledgements}
The authors wish to thank Giuseppe Lapadula for fruitful discussions on the
subject of this paper.

\end{document}